\documentclass[prl,amsmath,amssymb,floatfix,superscriptaddress,showpacs,footinbib]{revtex4}
\usepackage{graphicx}
\usepackage{dcolumn}
\usepackage{bm}
\usepackage{latexsym}
\usepackage{latexsym,graphics}
\usepackage{amscd} 
\usepackage[bookmarks, colorlinks=true, plainpages = false, citecolor = green, urlcolor = blue, filecolor = blue]{hyperref} 

\newcommand{\be}{\begin{equation}}
\newcommand{\ee}{\end{equation}}
\newcommand{\dif}{\,\mathrm{d}}

\DeclareMathOperator*{\argmax}{arg\,max}

\begin{document}

\title{Sampling with Costs}
\author{Joseph D. Skufca}
\email{jskufca@clarkson.edu}
\affiliation{Department of Mathematics, Clarkson University, Potsdam New York}%
\author{Daniel ben-Avraham}
\email{qd00@clarkson.edu}
\affiliation{Department of Physics, Clarkson University, Potsdam New York}%
\affiliation{Department of Mathematics, Clarkson University, Potsdam New York}%

\date{\today}

\begin{abstract}
We consider the problem of choosing the best of $n$ samples, out of a large random pool, when the sampling of each member is associated with a certain cost.  The quality (worth) of the best sample clearly increases with $n$, but so do the sampling costs, and one important question is how many to sample for optimal gain (worth minus costs).
If, in addition, the assessment of worth for each sample is associated with some ``measurement error," the perceived best out of $n$ might not be the actual best, complicating the issue.  Situations like this are typical in mate selection, job hiring, and food foraging, to name just a few.  We tackle the problem by standard order statistics, yielding suggestions for optimal strategies, as well as some unexpected insights.  
 
\end{abstract} 

\pacs{2.50.Le}

\maketitle

As a motivating example, consider the problem of the academic hiring committee when conducting a candidate search.  A large number of candidates apply, and, after filtering only to the highly qualified candidates, their  application records provide little insight into the multitude of issues that would determine which of these candidates is the ``best'' for the job, which is why we invite several of those candidates for a campus interview.  This process requires money, time, and effort, so (of course) we don't invite all candidates.  But how many should we invite?  Given that the evaluation process (the interview) does not provide perfect information about the eventual success of a candidate, how much good are we gaining by the interview?  If our first candidate does very well, should we make an offer, or should we wait to sample more from the pool?  If our initial slate of candidates was ``just Okay,'' what should we expect to gain by asking the Dean to let us invite more candidates?  The general difficulty is that one would like to choose ``the best,'' but if the candidates aren't very different, or if our ability to distinguish ``the best'' is not very good, then we may be wasting our resources.  


\section{Background}
We consider the problem of maximizing gain, on choosing an item from a large population pool.  First, imagine that we are presented with $n$ items, randomly selected from some population, where we would like to choose the item with the greatest {\it worth}, as measured by the value of some attribute which we denote as $A.$  We may treat the value of this attribute as a random variable, with distribution determined by the underlying population distribution, and denote the attribute value for the $i$th item as $A_i.$ Then 
\[
A_{\max}(n) = \max_{i=1,\ldots,n} A_i
\]
 is also a random variable, and the standard tools from {\it order statistics} may be applied to find the probability distribution for $A_{\max}$ in terms of the cumulative distribution function $F(a)$ of that attribute for the population.
Imagine further, that each {\it measurement} --- an assessment of the value of an item --- carries some {\it cost}, monetary or otherwise, so that the total cost of the measurements is $C_n$.  Then, the total {\it gain} to be gotten from the process is $g(n)=A_{\max}(n) -C_n$.  One important goal is to find $n$ which maximizes $g(n)$: how many people should one interview before hiring, how many mates should we date before proposing, how many cars to test-drive before buying? etc.  This problem is treated in Section~{\bf x}.

A complication arises when the evaluation process of the worth of each item is imperfect, yielding a somewhat erroneous value.  In that case, the perceived ``best" item out of $n$ might not coincide with the actual best, and that diminishes the expected gain.  The precise effect of noisy measurement, and how to work out an optimal strategy despite it, is treated in Section~{\bf y}.   The ubiquitous case where the items' worth and the error in measurement are each normally distributed is particularly enlightening, yielding some simple closed-form formulas, and we use it to demonstrate the general procedure.

Some further insights are developed in Section~{\bf z}, where we show that it always pays to sample three items, if it pays to sample at all, when the worth distribution and error distribution are both normal.
We conclude and discuss our findings in Section~{\bf w}.

\section{Expected gain with ideal measurement}
\subsection{Order Statistics and Worth}
We begin with the ideal case that the value of each item is assessed perfectly, without any measurement error.
Consider then a sample of $n$ i.i.d.~random variables $X_i$, taken from the distribution $p(x)$ --- the probability density function for the worth of our items --- which are reordered according to their ascending worth: $X_{(1)},X_{(2)},\dots,X_{(n)}$. Standard order statistics gives us the cumulative distribution function (cdf) for $X_{(k)}$:
\be
\label{eqstd_order} 
\Psi_{X_{(k)}}=\mathbb{P}(X_{(k)} \leq x)=\sum_{j=k}^n {n \choose j} P(x)^j (1-P(x))^{n-j},
\ee
where $P(x)=\int_{-\infty}^xp(x')\dif x'$ is the cdf of the items' worth.  Focusing on the largest item selected, we have
\be
\label{eqstd_largest} 
\Psi_{X_{(n)}}=\mathbb{P}(X_{(n)} \leq x)= P(x)^n\,,\qquad \psi_{X_{(n)}}(x)=nP(x)^{n-1}p(x)\,,
\ee
where the probability density function (pdf)  $\psi_{X_{(n)}}(x)$ was computed by differentiation.  A quick, alternative way to obtain this last result is by realizing that $P(x)$ denotes the probability that any of the $X_i$ be smaller than $x$.  Then, for the maximal value to be $x$, we need one of the $X_i$ to equal $x$, say $X_m=x$, while $X_j<x$ for $j\neq m$.  This happens with probability $nP(x)^{n-1}p(x)$, since $m$ can be chosen in $n$ different ways.  The expected value for this maximal order statistic, which we denote as $K_n,$ is computed as
\be 
\label{eqstd_exp}
K_n :=\mathbb{E}[X_{(n)}] = n \int_{-\infty}^{\infty} x P(x)^{n-1} p(x) \dif x\,.
\ee
A simple variable transformation shows that the analogous result for $p(x)'=ap(ax+b)$, is $K_n'=\frac{1}{a}(K_n-b)$.

For the flat distribution: $p(x)=1$ for $0<x<1$ (and zero otherwise), for example, one obtains $K_n=n/(n+1)$.
In general, however, no closed form solution exists for $K_n$, but numerical approximations for some  distributions can be found in most texts on order statistics and are available in statistical software packages.   For example, in the special case of standard normal variables, these expectations are called {\it rankits}, with these values required to make Q-Q plots.  

For large $n$, a simple, useful approximation for $K_n$, due to Van der Waerden, is given by $\int_{-\infty}^{K_n}p(x)\dif x\approx n/(n+1)$. (It does give the exact result for the flat distribution of the example.)  For the normal distribution, $\phi(x)=\frac{1}{\sqrt{2\pi}}e^{-x^2/2}$, this approximation yields $K_n\sim2\sqrt{\ln n}$. The very slow increase of $K_n$ with $n$ is quite typical, with the exception of fat-tailed distributions: for $p(x)=\alpha x^{-1-\alpha}$, $x>1$ (and zero elsewhere), for example, $K_n\sim n^{1/\alpha}$, which increases rapidly for small values of $\alpha$.  

\subsection{Costs, Gains, and Optimization}
The value of  $K_n$ is an increasing function of $n$, so if the goal is to ``get the very best," the strategy is simply to sample as many as possible.  In practical situations, however, there is invariably a {\it cost} associated with the sampling and measuring process: Bringing in candidates for interviews costs money and time; in the animal kingdom, courting many potential partners costs energy and time, delaying an eventual union and diminishing the chances for reproduction; or searching for the larger fruits exposes a forager to increasing danger, the longer the search, etc.  As a decision problem, the choice of $n$ should be based on what gives the most {net} benefit, or {\it gain}.  

Denote the { cost} of measuring the $i$th item by $c_i$, with cumulative cost
\begin{equation}
C_n = \sum_{i=1}^n c_i\,.
\end{equation}
The optimal sample size $n^*$ would then be given by the optimization problem
\[
n^* = \argmax_n \left(K_n -C_n \right)\,.
\]
In general, a reasonable assumption might be that the total cost is proportional to the number of samples, with fixed cost $c$ per item, i.e., $c_i=c$ and $C_n=nc$. We shall proceed under this assumption.  Note, however, that for the case of $n=1$, when only one item is picked, there is no point in measurement, since by necessity that one item is the best available.  Hence, we must also stipulate that $C_1=0$ (rather than $c$).

\begin{figure}[htb]
\includegraphics[width=.49\textwidth]{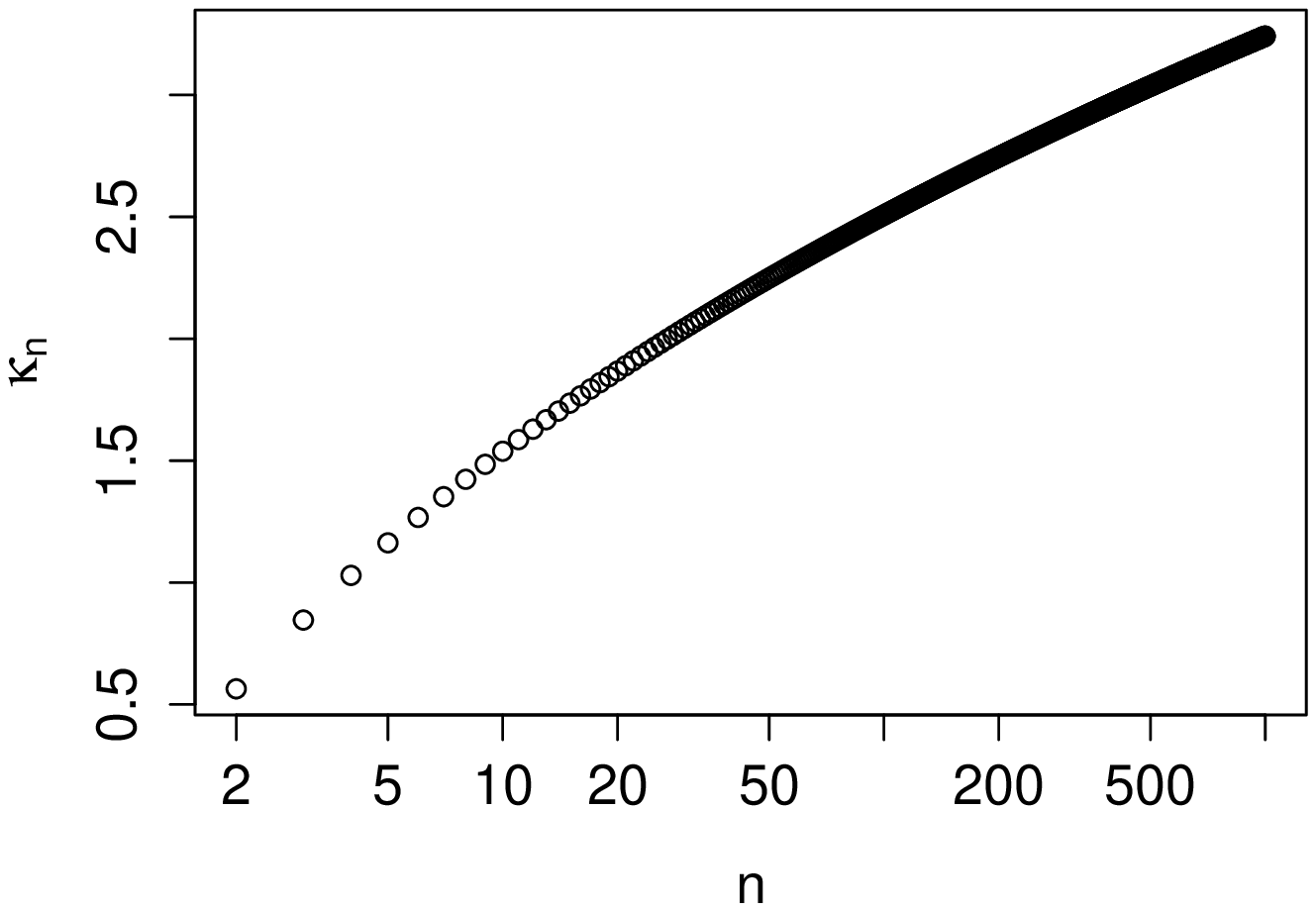} 
\includegraphics[width=.49\textwidth]{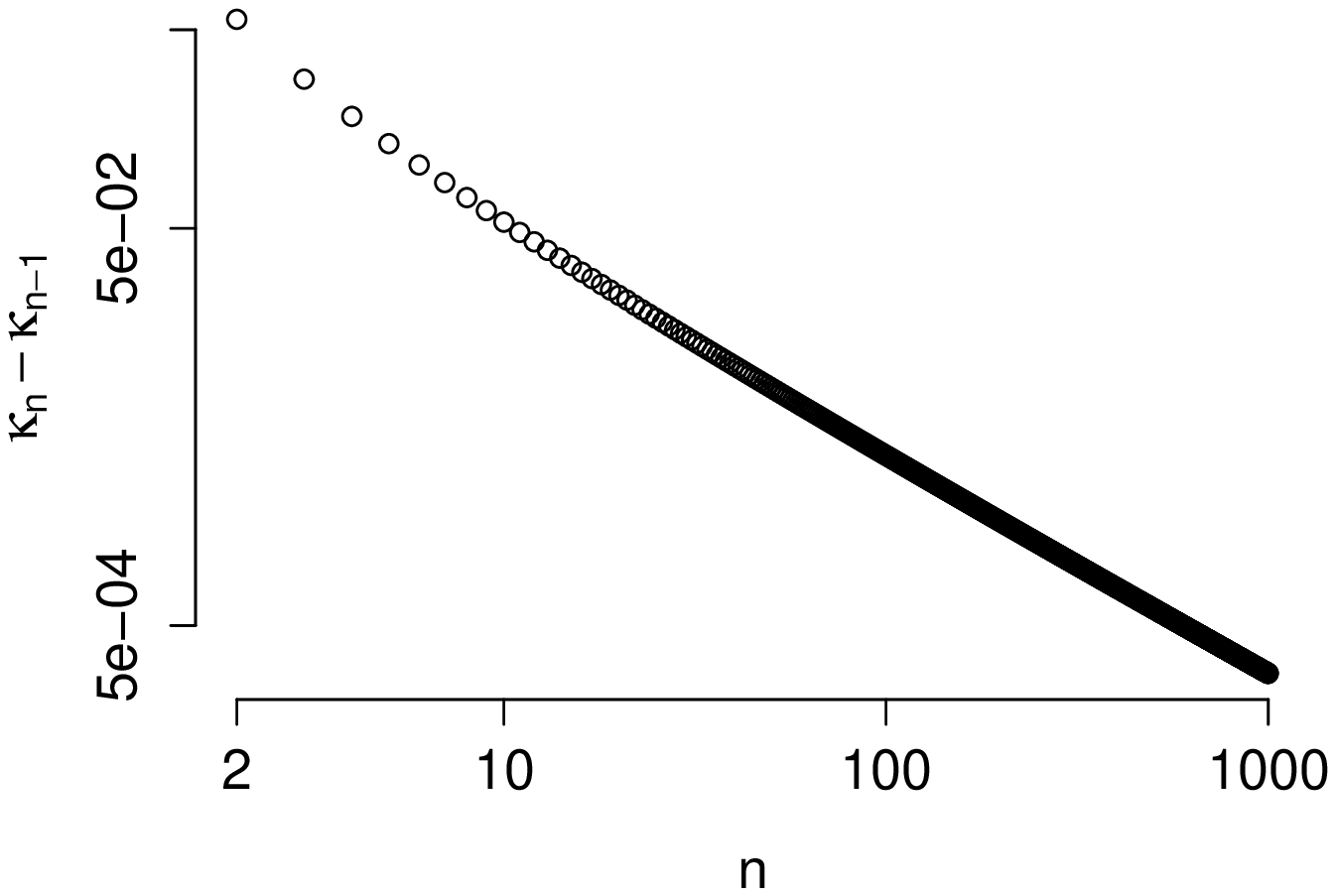}
\caption{{\bf Standardized Gain.}  Plots of $\kappa_n,$ the expected maximal statistic for the standard normal (Left),  and a log-scale plot of marginal $\kappa_n-\kappa_{n-1}$ (Right)  show that after a ``few'' samples, the gain grows very slowly with respect to the number of items examined.\label{figKn}}
\end{figure}

The marginal worth for sampling item $n$, given by
\[
k_n:=K_n-K_{n-1}\,,
\]
is usually a decreasing function in $n$, as illustrated by Figure~\ref{figKn}, for the normal distribution.  The optimal sample size  $n^*$ is then chosen as the largest $n$ such that the marginal worth exceeds the marginal cost:
\begin{equation}
\label{eq:nstarmarginal}
k_{n^*} > c \geq k_{n^* +1}\,,\qquad n^*\geq2\,.
\end{equation}
If $k_2<c$, the best strategy is to pick {\it one} item at random and keep it, without bothering to measure, as already discussed above.

As a simple example, consider $p(x)=1/a$ for $0<x<a$ (and zero elsewhere), for which $k_n=\frac{a}{n(n+1)}$.
If, furthermore, $a\gg c$, then $n^*\approx\sqrt{a/c}$.  For the normal distribution, $\phi_a(x)=\frac{1}{\sqrt{2\pi a^2}}e^{-x^2/2a^2}$, we get from Van der Waerden's approximation, $k_n\approx\frac{a}{n\sqrt{\ln n}}$ (for $n\gg1$), so
$n^*\approx a/c$, if $a\gg c$.  Finally, for the  freak case of a fat-tailed distribution, such as $p(x)=\alpha x^{-1-\alpha}$, $x>1$ (and zero elsewhere), $k_n\sim n^{(1/\alpha)-1}$, so for $\alpha<1$ the gain increases {\it indefinitely} with $n$, regardless of the mounting costs.

\section{The effect of measurement error \label{sec_err}}
We now turn to the case when the measurement of each item is not perfect, but associated with some error.
For simplicity and concreteness, throughout the remainder of the paper we focus on the most common scenario, where the worth of the items and
the error made in each measurement can both be described by the normal distribution.  The general case can be treated in much the same way, but is less transparent, since it is then impossible to push the analytical calculations as far.  

 To avoid any confusion, we denote the normal distribution of zero mean and variance $\sigma$ by $\phi_\sigma(x)$, instead of $p(x)$.  For $\sigma=1$, we simply use $\phi(x)$, dropping the subscript.
Likewise, we denote the expected maximal statistics of $\phi$ by $\kappa_n$ (instead of $K_n$).  Note that 
the expected maximal statistics for $\phi_{\sigma}$ is then $\sigma\kappa_n$.

Assume then that the $A_i$'s are independent and normally distributed with mean $\mu$ and variance $a$. We define the {\it return} for the $i$th item as
\be 
X_i = A_i-\mu\,,
\ee
where $\mu$ is the worth mean. The $X_i$'s are then i.i.d.\ random variables, described by the normal pdf
$\phi_a(x)=\frac{1}{\sqrt{2\pi a^2}}e^{-x^2/2a^2}$.  Assume, further, that each measurement is associated with an error
$Y_i$, and that the $Y_i$'s are i.i.d.\ random variables described by the normal pdf $\phi_b(y)$.  Thus, the actual value {\it measured} for the $i$th item is
\be 
W_i := X_i + Y_i\,.
\ee
From standard results, we observe that $Z_i$ is normally distributed, with mean $0$ and variance $\sqrt{a^2+b^2}$   (see eq.~(\ref{eq_Z}), below).

Because the process is independent for each $i,$ we (for the moment) drop the subscript in order to (notationally) ease the discussion in understanding the relationship between the {measured value} $W$ and the actual {return} $X$.
We first ask {\it what is the expected return given a particular measured value?}  Because sample and measurement error are independent, the joint distribution for $X$ and $Y$ is given by
\be
f_{XY}(x,y)=\phi_a(x)\phi_b(y)=\frac{1}{{2\pi ab}} \; e^{ -(\frac{x^2}{2a^2} +\frac{y^2}{2b^2}) }\,.
\ee
We perform a change of variables, $Y=W-X$ to find the joint distribution of $X$ and $W,$
\be 
\label{eqjoint_XW}
f_{XW}(x,w)=\phi_a(x)\phi_b(w-x)=\frac{1}{{2\pi ab}} \; e^{ -\frac{x^2}{2a^2} -\frac{(w-x)^2}{2b^2} }\,,
\ee
remarking that the substitution is simplified by the fact that $\partial Y / \partial W= 1$.
The pdf for measured values $W$ is then
\be
\label{eq_Z}
f_W(w)=\int_{-\infty}^{\infty}f_{XW}(x,w)\dif x=\frac{1}{\sqrt{2\pi(a^2+b^2)}}\;e^{-\frac{w^2}{2(a^2+b^2)}}=\phi_{\sqrt{a^2+b^2}}(w)\,,
\ee
while the
pdf for a return $x$, conditioned on a measured value $w$, is then
\be 
f_X (x|W=w) = \frac{f_{XW}(x,w)}{f_W(w)}=\frac{1}{\sqrt{2\pi\sigma^2}} \; e^{ -\frac{\left(x- \frac{a^2}{a^2+b^2} w  \right)^2}{2\sigma^2} }\,,
\ee
where $\sigma =\frac{ab}{ \sqrt{a^2+b^2}}.$
The required conditional expectation is then easily obtained:
\be 
\label{eqexp_cond_x}
\mathbb{E}[X|W=w] = \int_{-\infty}^{\infty} x f_X (x|W=w) \dif x = \frac{a^2}{a^2+b^2}\, w:=\eta^2 w\,.
\ee

Armed with this result we can now complete our original goal, of determining the expected return on selecting the {largest} item based on {\it measured} values.  Formally, we define
\[
Q=X_k
\]
where $k$ satisfies
\[
W_k =W_{(n)}\,,
\]
the largest order statistic of the sample measured values.  Stated more directly, $Q$ represents the return of the item that {\it measured} to be the largest.
We may compute the required expectation of $Q$ by integrating (\ref{eqexp_cond_x}) against $\psi_{W_{(n)}}$, the pdf for $W_{(n)}$ --- obtained from (\ref{eq_Z}) and (\ref{eqstd_largest}) --- to yield,
\be
\label{eq:EQ1}
\begin{split}
\mathbb{E}[Q] & = \int_{-\infty}^{\infty} \mathbb{E}[X|W=w]  \psi_{W_{(n)}}(w) \dif w \\
& = \frac{a^2}{a^2+b^2} \mathbb{E}[W_{(n)}] \\ 
&= \frac{a^2}{a^2+b^2} \sqrt{a^2+b^2} \kappa_n  \\ 
& =(a\kappa_n)\frac{a}{\sqrt{a^2+b^2}}=\eta (a\kappa_n)\,.
\end{split}
\ee
Note that $a\kappa_n$ is the result one expects in the ideal case, when there is no measurement error.
The net effect of measurement error, then, is to degrade the gain that could be obtained in the ideal case, by the factor $\eta=a/\sqrt{a^2+b^2}$.  (This would not be the case for distributions other than normal, in general, but one expects qualitatively similar behavior.) For $b\ll a$, $\eta\approx 1-\frac{b^2}{2a^2}$ and the degradation is minimal (and vanishing as $b\to 0$).  For $b\gg a$, however, $\eta\approx a/b$ and the degradation is large.  The latter case explains our  ``common sense'' understanding of two common situations:
\begin{itemize}
\item If there is not much  difference between items ($a$ small), don't bother to measure, just pick one.
\item If you {\it can't  tell} the difference between items ($b$ large),  don't bother to measure, just pick one.
\end{itemize}

\section{So, how many should we try?}
Well, we have already  answered this question, formally, by providing a way to compute $n^*$, the optimal number of trials.  But often one's search is less well planned, or the optimal strategy cannot be followed, due to external constraints (e.g., the funding for bringing interviewees on campus comes from your Dean).  Here 
we develop two important strategies to help deal with such problems.  The first strategy ignores pre-planning, and addresses the immediate question whether to sample once more, based on what we already have at hand.  The second strategy establishes a reasonable minimum of tries when one is pressed to terminate the searching prematurely.  

\subsection{Should we try one more?}
The analysis leading to the criterion of (\ref{eq:nstarmarginal}) addresses the question of how many items to sample, based on careful and deliberate {\it a priori} planning.  
In many instances, however, the sampling process is not pre-planned, but {\it sequential} (e.g., should I try on one more pair of jeans before making my purchase).  In such instances, the decision whether to sample one more is based only on the current information --- the measured value of the current best choice.  

Suppose that after some amount of sampling our best choice measures to be $w_0.$  If we sample one additional item, with measured value $W$, then we would prefer the old sample if $W \leq w_0,$ but switch to the new if $W>w_0.$  Then the expected increase in worth, conditioned on $W=w$, would be given by
\begin{equation}
h(w)=
\begin{cases}
0 &  \mbox{  if  } w \leq w_0, \\
\eta^2 (w-w_0)  & \mbox{  for  } w>w_0,
\end{cases}
\label{eq:condgain}
\end{equation}
where we have applied the result of \eqref{eqexp_cond_x}.   The unconditional expectation of gain on sampling one more, $V^+_{w_0}$, is then given by 
\begin{equation}
V^+_{w_0}=\int_{-\infty}^{\infty} h(w) f_W(w) \dif w=\int_{-\infty}^{\infty} h(w) \phi_{\sqrt{a^2+b^2}}(w) \dif w\,.
\end{equation}
Carrying out the integrals, and expressing the final result in terms of the standard normal distribution (with unit
variance), we obtain
\be
V^+_{w_0}=a\eta[\phi(z_0)+z_0(\Phi(z_0)-1)]:=a\eta v^+(z_0)\,,
\ee
where $z_0=w_0/\sqrt{a^2+b^2}$.    As a general guideline, $v^+(z)$ is a rapidly decreasing function of $z$:  $v^+(z)\sim |z|$, for $z\ll0$, $v^+(0)=1/\sqrt{2\pi}$, and $v^+(z)\sim e^{-z^2/2}$, for $z\gg0$. The (sequential) decision whether to sample one more is  based on whether $V_{w_0}^+>c$ (sample!) or not.  In Figure~xx, we plot the function $v^+(z)$ used for making this decision.

\subsection{Try at least three, or none!}
Under some circumstances, there is external pressure to limit the sampling to a small number of items, sometimes even when it is clear that a longer search would be more advantageous. For example, the funding for the search, such as in the case of hiring new faculty, might come from an external source (the Dean) , and one faces pressure to terminate the process as early as possible.  We here answer the question ``What is a reasonable minimum amount of tries?"  relevant to such situations.

We assume that the pdf of the items' value is normal, with average $\mu$ and variance $a$, and that the pdf of the error in measurement is also normal, with variance $b$ (and zero average).  Assume furthermore that the cost of measuring each item is $c$.  Then, if $\mu<c$, it pays to simply pick one item, at random, without measuring.  The expected gain in that case is $g_1=\mu$.  

Does it pay, instead, to try two items?  According to our results for selecting with measuring errors, the expected maximal worth of two items is $a\eta\kappa_2+\mu$, so that the expected gain is $g_2=a\eta\kappa_2+\mu-2c$.  Thus, it pays to try two items if $g_2>g_1$, or $a\eta\kappa_2-2c>0$.

We shall now prove that $\kappa_3=\frac{3}{2}\kappa_2$.  In that case, the expected gain from trying three items at the outset is
$g_3=a\eta\kappa_3+\mu-3c=g_2+(a\eta\kappa_2-2c)/2$.  Thus, whenever it pays to try two items, it does pay
even more to try three!  This suggest the following ``minimalist" strategy: {\it If you believe that the cost of measuring is too high for even a small number of items, then just pick one at random (and don't bother to measure). Otherwise,  try at least three}.

Using the result~(\ref{eqstd_exp}), and exploiting the fact that $\phi(x)$ is an even function of $x$, while $x$ and $\Phi(x)-\frac{1}{2}$ are odd, the proof is straightforward:
\[
\begin{split}
\kappa_3 &= 3\int_{-\infty}^{\infty}x\,\Phi(x)^2\phi(x)\dif x=3\int_{-\infty}^{\infty}x\left[\left(\Phi(x)-\frac{1}{2}\right)+\frac{1}{2}\right]^2\phi(x)\dif x\\
&=3\int_{-\infty}^{\infty}x\left[2\left(\Phi(x)-\frac{1}{2}\right)\cdot\frac{1}{2}\right]\phi(x)\dif x=3\int_{-\infty}^{\infty}x\,\Phi(x)\phi(x)\dif x\\
&=\frac{3}{2}\kappa_2\,.
\end{split}
\]
Incidentally, the above proof also shows that $K_3=\frac{3}{2}K_2$ for any pdf that is an even function of its argument,
and  the same symmetry trick can be used to  obtain $K_{2n+1}$ in terms of $K_2,K_4,\dots,K_{2n}$; for example, $K_5=\frac{5}{2}K_4-\frac{5}{2}K_2$, etc.

\section{Discussion and Conclusion}

As primary results from this paper, we briefly restate what we consider as the key analytic contributions:
\begin{enumerate}
\item If we are measuring with error, and determine a particular measured value $w,$ then the expected {\it true value} (accounting for stochastic differences in the population, not error in our measurement) if given be 
\be 
\mathbb{E}[X|W=w] = \int_{-\infty}^{\infty} x f_X (x|W=w) \dif x = \frac{a^2}{a^2+b^2}\, w:=\eta^2 w\,.
\ee
\item If we intend to measure $n$ items and select the item that measures as the best, the expected benefit of that process is given by
\begin{equation}
V(n,a,b):=\kappa_n \frac{a^2}{\sqrt{a^2+b^2}}= \frac{a \kappa_n}{\sqrt{1+(b/a)^2}}=\eta a\kappa_n \,.
\end{equation}

\item If we currently have an item which measures $w_0,$ then the expected gain in worth, on picking one more item to measure, is given by
\be
V^+_{w_0}=a\eta[\phi( z_0)+ z_0(\Phi( z_0)-1)]:=a\eta v^+( z_0)\,,
\ee
where $ z_0=w_0/\sqrt{a^2+b^2}$.
\end{enumerate}

In the context of our original motivating example (the candidate search), we note that item~(2) addresses the question of how many people the Dean might let us invite, but that decision would still require some means of determining {\it costs} of a candidate visit measured in the same units as the {\it value} of selecting a better candidate.  Item~(3) addresses the question of whether we should make an offer to our current ``best candidate,'' or should we wait to see another candidate.  Item~(1) relates directly to the issue of the importance of having a good measuring system - the interview process itself, where we would remark that $m$ can be reduced through repeated measuring, equivalent to requiring the candidate to stay for a longer visit and conduct more interviews.  However, it is worthwhile to note that the relationship between the math and our illustrative problem is mostly qualitative, in that our normality assumptions, as well as the idea that we have some idea of mean and variance of the population and our measuring device is not reasonable.

As a component of discussion, we think it is worthwhile to comment upon the implications of these results.  We recall that the sampling process can be assumed to have costs, so decision theory principles drive lead us to the simple conclusion that {\it we should only sample more items if the expected gain is less than the cost of sampling.}  Consequently, our analytic formulas provide additional insight into the process.
\begin{itemize}
\item  The more we sample, the better should be our performance in selection, so long as we do not exceed where sampling costs exceed expected benefits.
\item  If our measurement system is not very accurate, we suffer two effects.  On the one hand, we are less able to select the best item, but, additionally, are expected gain is reduced.  For a fixed marginal cost to sample, that means we will stop sampling sooner, settling earlier in the process, further reducing are likelihood of finding an ``exceptionally good'' item.
\item As corollary, if we want to find very good items, sampling costs must be very low.
\item As second corollary, if we can reduce are measurement error, it can become cost effective to sample more items.
As a numerical example, if sampling cost was such that we would have looked at $n=10$ items, are standardized expected gain is $K_{10} \approx 1.54.$ If the per item sampling cost were reduced by a factor of 10, then based on the marginal benefit being greater than marginal cost, we would sample $n=70$ items, with $K_{70}\approx2.38.$  Based on the rapid decay, we would note that   
the benefit grows roughly with like $\log n.$
\end{itemize}

If we examine these principles playing out in arenas such as mate selection, we would (perhaps) have to ignore the competitive aspect (your proposed mate must also choose you over other possible mate choices.  However, one could use these results to infer that {\it performance in the mate selection arena is enhanced if ``dating'' is cheap.}  Specifically, if we want to find a very good mate, then we must follow be willing to perform more sampling.  Biologically, there is an inherent risk cost associated with moving from one mate choice to another.  We note that there appears to have been evolutionary pressure in this direction \cite{buss2003} as it is a well observed phenomena that the body reacts (hormonally) to provide increased pleasure during the first stages of a relationship (the thrill of dating).  Perhaps this pleasure boost should be viewed as {\it decreasing} the cost associated with sampling so that there is marginal reason to sample additional items before choosing a mate.

•


\end{document}